\documentclass{article}
\usepackage{latexsym,amsfonts,amssymb}
\usepackage[latin1]{inputenc}
\usepackage{amsmath}
\usepackage{enumitem}

\title{\sc Multiplying a conjugacy class by its inverse in a finite group}

\author{Antonio Beltr\'an\\
\footnotesize
Departamento de Matem\'aticas,\\
\footnotesize Universidad Jaume I, \footnotesize
12071 Castell\'on, Spain\\
\footnotesize
e-mail: abeltran@mat.uji.es\\
\\
Mar\'{\i}a Jos\'e Felipe\\
\footnotesize
Instituto Universitario de Matem\'atica Pura y Aplicada,\\
\footnotesize Universidad Polit\'ecnica de Valencia, \footnotesize
46022 Valencia, Spain\\
\footnotesize e-mail: mfelipe@mat.upv.es \\
\\ Carmen Melchor\\
\footnotesize
Departamento de Matem\'aticas,\\
\footnotesize Universidad Jaume I, \footnotesize
12071 Castell\'on, Spain\\
\footnotesize
e-mail: cmelchor@uji.es      }

\date{}
\begin{document} \maketitle

\begin{abstract}
Suppose that $G$ is a finite group and $K$ a non-trivial conjugacy class of $G$ such that $KK^{-1}=1\cup D\cup D^{-1}$ with $D$ a conjugacy class of $G$. We prove that $G$ is not a non-abelian simple group. We also give arithmetical conditions on the class sizes determining the structure of $\langle K\rangle$ and $\langle D\rangle$. Furthermore, if $D=K$ is a non-real class, then $\langle K\rangle$ is $p$-elementary abelian for some odd prime $p$.



\end{abstract}

\bigskip

\section{Introduction}
There are many studies about the structure of a finite group focused on the product of its conjugacy classes. Perhaps, the most relevant problem was posed by Z. Arad and M. Hergoz (\cite{products}) who conjectured that if $S$ is a non-abelian simple group and $A$ and $B$ are non-trivial conjugacy classes of $S$, then $AB$ (defined as the set $\lbrace ab\,  | \, a\, \in\,  A, \,  b\,  \in\,  B\rbrace$) cannot be a single conjugacy class of $S$. This conjecture still remains open although some specific cases have been solved. For instance, in \cite{MooriViet}, the conjecture is verified for several families of finite simple groups of Lie type. On the other hand, Arad and E. Fisman proved in \cite{AradFisman} that if $C$ and $D$ are non-trivial conjugacy classes of a finite group $G$ such that either $CD=C\cup D$ or $CD=C^{-1}\cup D$, then $G$ is not a simple group. However, no solvability information of the subgroups $\langle C\rangle$ or $\langle D\rangle$ was given. Recently, G. Navarro and R.M. Guralnick (\cite{GuralnickNavarro}) have proved that when a conjugacy class $K$ of a finite group $G$ satisfies that $K^2$ is a conjugacy class, then $\langle K\rangle$ is a solvable (normal) subgroup of $G$ by appealing to the Classification of the Finite Simple Groups (CFSG). This is, of course, consistent with Arad and Hergoz's conjecture. \\

In general, for any $G$-invariant subset $X$ of $G$, we denote by $\eta(X)$ the number of distinct conjugacy classes appearing in $X$. Let $K$ be a conjugacy class of a finite group $G$. When we multiply $K$ by its inverse class, $K^{-1}$, then $KK^{-1}$ is a $G$-invariant set. We will prove that if $\eta(KK^{-1})=2$, then $G$ is not simple. The fact that $\eta(KK^{-1})=3$ does not imply that $\langle K\rangle$ or $\langle KK^{-1}\rangle$ is solvable. Really, $\langle KK^{-1}\rangle$ may be even simple. For instance, if $G=S_n$ for any $n\geq 5$ and $K$ is the conjugacy class of transpositions, then $\eta(KK^{-1})=3$ and $\langle KK^{-1}\rangle=A_{n}$. In this paper we study the particular case in which $KK^{-1}=1\cup D\cup D^{-1}$ with $D$ a conjugacy class of $G$, and we demonstrate that $G$ cannot be simple by means of the CFSG and a character theoretical property characterizing such condition for conjugacy classes.\\

\textbf{Theorem A.} {\it Let $K$ be a non-trivial conjugacy class of a finite group $G$ and suppose that $KK^{-1}= 1 \cup D\cup D^{-1}$, where $D$ is a conjugacy class of $G$. Then $G$ is not a non-abelian simple group. In particular, this theorem holds if $KK^{-1}= 1 \cup D$.}\\

We remark that $KK^{-1}=1 \cup D$ forces that $D$ is real but it does not necessarily imply that $K$ is a real class too (see section 5). Moreover, under the assumption of Theorem A, if $K$ is real, then we see (Lemma 3.1) that $D$ is real too and thus, $K^{2}=1\cup D$. In this case, the structure and solvability of $\langle K\rangle$ is obtained in \cite{Nuestro}, without employing the CFSG.\\

In order to prove Theorem A, we use the following characterization in terms of characters of the property appearing in such theorem. We denote by Irr$(G)$ the set of all irreducible complex characters of $G$.\\

\textbf{Theorem B.} {\it Let $G$ be a group and $x, d\in G$. Let $K=x^G$ and $D=d^G$. The following are equivalent:}
\begin{enumerate}[label=\alph*)]
\item $KK^{-1}=1\cup D \cup D^{-1}$
\item For every $\chi \in$ {\rm Irr}$(G)$ $$|K||\chi(x)|^2=\chi(1)^2+\frac{(|K|-1)}{2}\chi(1)(\chi(d)+\chi(d^{-1})).$$
\end{enumerate}

{\it In particular, if $D=D^{-1}$, then $KK^{-1}=1\cup D$ if and only if for every $\chi \in$ {\rm Irr}$(G)$ $$|K||\chi(x)|^2=\chi(1)^2+(|K|-1)\chi(1)\chi(d).$$}

Under the hypotheses of the particular case of Theorem A, we have only been able to prove the solvability of the subgroup $\langle K\rangle$ in some concrete cases.\\

\textbf{Theorem  C.} {\it Let $K$ be a conjugacy class of a finite group $G$ and suppose that $KK^{-1}= 1 \cup D$, where $D$ is a conjugacy class of $G$. Then $|D|$ divides $|K|(|K|-1)$ and $\langle K\rangle/\langle D\rangle$ is cyclic. In addition,
\begin{enumerate}
\item If $|D|=|K|-1$, then $\langle K\rangle$ is metabelian. More precisely, $\langle D\rangle$ is $p$-elementary abelian for some prime $p$.
\item If $|D|=|K|$, then $\langle K\rangle$ is solvable with derived length at most 3.
\item If $|D|=|K|(|K|-1)$, then $\langle K\rangle$ is abelian.
\end{enumerate}}

We will provide examples showing that each case is feasible as well as an example in which $|D|$ is a divisor of $|K|(|K|-1)$ distinct from those appearing in Theorem C and satisfying $dl(\langle K\rangle)=3$. Although in these examples $\langle K\rangle$ is solvable, the general proof remains open.\\

Observe that if $KK^{-1}=1\cup K$, then $K$ is real, so, $K^{2}=1\cup K$, and it easily turns out that $\langle K\rangle$ is $p$-elementary abelian since $\langle K\rangle$ is a minimal normal subgroup and all the non-trivial elements have the same order. Thus, we wonder what happens when $KK^{-1}=1\cup K\cup K^{-1}$ and $K\neq K^{-1}$. The interesting point of this non-real case is that more arithmetical information can be obtained by simply working within the complex group algebra $\mathbb{C}[G]$.\\

\textbf{Theorem D.} {\it Let $K$ a non-real conjugacy class of $G$ satisfying $KK^{-1}=1\cup K\cup K^{-1}$. Then $\langle K\rangle $ is an elementary abelian $p$-group with $p$ an odd prime such that if $|\langle K\rangle|=p^n$, then $p\equiv 3$ $mod$ $4$ and $n$ is odd. Moreover, $|K|=\frac{p^n-1}{2}$.}

\section{Preliminary results and proof of Theorem B}
We begin this section by presenting several results appeared in the literature that we need in order to prove Theorem A. The next theorem will be useful in the interest of discarding the alternating groups $A_{n}$ with $n> 5$ in Theorem A.\\

\textbf{Theorem 2.1} (Theorem A of \cite{AdanBante1}) {\it Let $S_n$ be the symmetric group of $n$-letters, $n>5$ and $\alpha$, $\beta$ $\in$ $S_{n}\setminus 1$. Then $\eta(\alpha^{S_{n}} \beta^{S_{n}})\geq 2$ and, if $\eta(\alpha^{S_{n}}\beta^{S_{n}})=2$ then either $\alpha$ or $\beta$ is a fixed point free permutation. Assume that $\alpha$ is fixed point free. Then one of the following holds
\begin{enumerate}
\item $n$ is even, $\alpha$ is the product of $n/2$ disjoint transpositions and $\beta$ is either a transposition or a $3$-cycle.
\item $n$ is a multiple of $3$, $\alpha$ is the product of $n/3$ disjoint $3$-cycles and $\beta$ is a transposition.
\end{enumerate}}

The following result due to R.M. Guralnick and G.R. Robinson is an extension for odd primes of Glauberman's ${\rm Z}^{*}$-Theorem \cite{Glauberman}.\\

\textbf{Theorem 2.2} (Theorem D of \cite{GurRob}) {\it Let $G$ be a finite group. If $x \in G$ has order $p$ and $[x, g]$ is a $p'$-element for every $g\in G$, then $x$ is central modulo {\rm \textbf{O}}$_{p'}(G)$. } \\

We give a variation of such theorem for $p$-elements (not necessarily of order $p$) by adding an hypothesis to the class size of the $p$-element. We use the following property (based on the CFSG) so as to obtain our variation, which is also a key result to prove Theorem 2.2.\\

\textbf{Lemma 2.3} (Theorem 4.1 of \cite{GurRob}) {\it Let $G$ be a finite group. If $x \in G$ has order $p$ and is not central modulo {\rm \textbf{O}}$_{p'}(G)$, then $x$ commutes with some conjugate $x^g\neq x$, for some $g\in G$. } \\

Our extension of Theorem 2.2 is the following.\\

\textbf{Lemma 2.4.} {\it Let $G$ be a finite group. Let $x \in G$ be a $p$-element such that $|x^G|$ is a $p'$-number and that  $[x, g]$ is a $p'$-element for every $g\in G$. Then $x$ is central modulo {\rm \textbf{O}}$_{p'}(G)$. }\\

 {\it Proof. } Let $o(x)=p^r$. By Theorem 2.2, we can assume that $r>1$. We write $y=x^{p^{r-1}}$, so $o(y)=p$. We will argue by induction on $|G|$.\\

We claim that there is no $g\in G$ such that $y\neq y^g\in {\rm \textbf{C}}_{G}(x)$. Suppose that there exists $g\in G$ such that $y^{g}$ centralizes $x$ and $y\neq y^g$. Since $|x^{G}|$ is a $p' $-number, we can choose $P\in$ Syl$_{p}(G)$ such that $P\subseteq {\rm \textbf{C}}_{G}(x) \subseteq {\rm \textbf{C}}_{G}(y)$, so there exists $n\in {\rm \textbf{C}}_{G}(y)$ such that $y^{gn}\in P \subseteq{\rm \textbf{C}}_{G}(x)$. In addition, $y^{-1}\in {\rm \textbf{C}}_{G}(x)$, so $[y, gn] \in {\rm \textbf{C}}_{G}(x)$, that is, $[gn, y, x]=1$. Moreover, $[y, x, gn]=1$, so by the Three Subgroups Lemma, we get $[x, gn, y]=1$. Then $x^{-1}x^{gn}=[x, gn]\in {\rm \textbf{C}}_{G}(y)$. As a result, $x^{gn}\in {\rm \textbf{C}}_{G}(y)$. Analogously, there is $n'\in {\rm \textbf{C}}_{G}(y)$ such that $x^{gnn'}\in {\rm \textbf{C}}_{G}(x)$ and hence, $x^{-1}x^{gnn'}$ is a $p$-element. By applying the hypothesis, $x=x^{gnn'}$, which implies that $gnn'\in {\rm \textbf{C}}_{G}(x)\subseteq {\rm\textbf{C}}_{G}(y)$. In particular, $g\in {\rm \textbf{C}}_{G}(y)$, a contradiction. As a consequence, there is no $g\in G$ such that $y\neq y^{g}\in {\rm \textbf{C}}_{G}(x)$ as claimed. So, $y{\rm \textbf{O}}_{p'}(G)\in {\rm \textbf{Z}}(G/{\rm \textbf{O}}_{p'}(G))$ by Lemma 2.3. \\

We distinguish two cases depending on whether ${\rm \textbf{O}}_{p'}(G)\neq 1$ or ${\rm \textbf{O}}_{p'}(G)=1$. Assume first that ${\rm \textbf{O}}_{p'}(G)=1$, and hence $y\in {\rm \textbf{Z}}(G)_{p}\neq 1$. Let $\overline{G}=G/{\rm \textbf{Z}}(G)_{p}$. We have that $1\neq \overline{x}$ is a $p$-element, $|\overline{x}^{\overline{G}}|$ is a $p'$-number and $[\overline{x}, \overline{g}]=\overline{[x,g]}$ is $p'$-element for every $\overline{g}\in \overline{G}$. By the inductive hypothesis, $\overline{x}{\rm \textbf{O}}_{p'}(\overline{G}) \in {\rm \textbf{Z}}(\overline{G}/{\rm \textbf{O}}_{p'}(\overline{G}))$. But observe that ${\rm \textbf{O}}_{p'}(\overline{G})=1$, so $ \overline{x}\in {\rm \textbf{Z}}(\overline{G})$. This means that $[x, g]\in {\rm \textbf{Z}}(G)_{p}$ for every $g\in G$. By hypothesis, $[x, g]$ is also a $p'$-element, so $[x ,g]=1$ for every $g\in G$. This shows that $x\in {\rm \textbf{Z}}(G)$ and the theorem is proved. \\

Suppose now that ${\rm \textbf{O}}_{p'}(G)\neq 1$ and let $\overline{G}=G/{\rm \textbf{O}}_{p'}(G)$. Again $1\neq \overline{x}$ is a $p$-element, $|\overline{x}^{\overline{G}}|$ is a $p'$-number and $[\overline{x}, \overline{g}]=\overline{[x,g]}$ is a $p'$-element for every $g\in G$. Notice that ${\rm \textbf{O}}_{p'}(\overline{G})=1$. By induction, $\overline{x} \in {\rm \textbf{Z}}(\overline{G})$, or equivalently $[x, g] \in {\rm \textbf{O}}_{p'}(G)$  for every $g\in G$, so the proof is finished. $\Box$\\

We state two very well-known Burnside's results.\\

{\bf Lemma 2.5} (Theorem 1.2.6 of \cite{Michler}) {\it A non-cyclic 2-group $P$ has only one involution if and only if $P$ is a generalized quaternion group.}\\

{\bf Lemma 2.6} (Lemma 15.1 of \cite{Huppert}) {\it Let $\chi \in$ {\rm Irr}$(G)$ and $K$ a conjugacy class of an element $g\in G$. Suppose that $(|K|, \chi(1))=1$. Then either $g\in {\rm \bf Z}(\chi)$, that is $|\chi(g)|=\chi(1)$, or $\chi(g)=0$.}\\

For our purposes we develop some properties related to the product of class sums in the complex group algebra $\mathbb{C}[G]$ of a finite group $G$ in order to prove Theorems B and D. Let $K_{1}, \ldots, K_{n}$ be the conjugacy classes of $G$. We denote by $\widehat{K_{i}}$ the class sum of the elements of $K_{i}$ in $\mathbb{C}[G]$. It is known that $\lbrace \widehat{K_{1}}, \ldots, \widehat{K_{n}}\rbrace$ is a $\mathbb{C}$-basis of ${\rm \textbf{Z}}(\mathbb{C}[G])$, the center of $\mathbb{C}[G]$. If $S$ is a $G$-invariant set, then $$\widehat{S}=\sum_{g\in S}g=\sum_{i=1}^{n}n_{i}\widehat{K_{i}}$$ denotes the sum of all elements in $S$ and we denote by $n_{i}=(\widehat{S}, \widehat{K_{i}})=(\widehat{K_{i}}, \widehat{S})$ the multiplicity of $\widehat{K_i}$ in $\widehat{S}$, which is of course a non-negative integer. For more details, we refer the reader to Chapter 3 of \cite{Huppert}.\\

\textbf{Lemma 2.7.} {\it If $D_{1}$, $D_{2}$ and $D_{3}$ are conjugacy classes of a finite group $G$, then
\begin{enumerate}
\item $(\widehat{D_{1}}\widehat{D_{2}}, \widehat{D_{3}})=(\widehat{D_{1}^{-1}}\widehat{D_{2}^{-1}}, \widehat{D_{3}^{-1}})$
\item $(\widehat{D_{1}}\widehat{D_{2}}, \widehat{D_{3}})=|D_{2}||D_{3}|^{-1}(\widehat{D_{1}}\widehat{D_{3}^{-1}}, \widehat{D_{2}^{-1}})$
\item $(\widehat{D_{1}}\widehat{D_{2}}, \widehat{D_{1}})=|D_{2}||D_{1}|^{-1}(\widehat{D_{1}}\widehat{D_{1}^{-1}}, \widehat{D_{2}^{-1}})=(\widehat{D_{2}}\widehat{D_{1}^{-1}}, \widehat{D_{1}^{-1}})=(\widehat{D_{2}^{-1}}\widehat{D_{1}}, \widehat{D_{1}})$.
\end{enumerate}}

{\it Proof.} See the proof of Theorem A of \cite{AradFisman}. $\Box$\\

{\it Proof of Theorem B.} Suppose that $KK^{-1}=1\cup D\cup D^{-1}$. Observe that by Lemma 2.7(1) $(\widehat{K}\widehat{K^{-1}}, \widehat{D})=(\widehat{K^{-1}}\widehat{K}, \widehat{D^{-1}})=(\widehat{K}\widehat{K^{-1}}, \widehat{D^{-1}})$. Thus, $\widehat{K}\widehat{K^{-1}}=|K|\widehat{1}+m\widehat{D}+m\widehat{D^{-1}}$ where $m$ is a positive integer. Therefore,
\begin{equation}\label{1}
|K|^{2}=|K|+2m|D|.
\end{equation} By applying Theorem 3.9 of \cite{Huppert},
$$\frac{|K|^{2}|\chi(x)|^{2}}{\chi(1)^{2}}=|K|+\frac{m|D|\chi(d)}{\chi(1)}+\frac{m|D|\chi(d^{-1})}{\chi(1)}$$
 for each $\chi \in$ Irr$(G)$. Taking into account Eq.(\ref{1}) and rearranging the equality we get the stated formula in b). Note that in the case $D=D^{-1}$, instead of Eq.(\ref{1}), we obtain $|K|^{2}=|K|+m|D|$ for some positive integer $m$ and by reasoning as before we have the result. \\

Conversely, suppose that b) is true. Let $C_{i}$ be with $1\leq i\leq n$ the conjugacy classes of $G$. By exercise 3.9 of \cite{Isaacs}, for any pair of conjugacy class sums $\widehat{C}_{m}$ and $\widehat{C}_{n}$ with representatives $c_{m}$ and $c_{n}$ we have
\begin{equation*}
\widehat{C}_{m}\widehat{C}_{n}=\sum_{k}\alpha_{k}\widehat{C}_{k}
\end{equation*} where
 \begin{equation*}
 \alpha_{k}=\frac{|C_{m}||C_{n}|}{|G|}\sum_{\chi\in {\rm Irr}(G)}\frac{\chi(c_{m})\chi(c_{n})\overline{\chi(c_{k})}}{\chi(1)}
\end{equation*} and $c_{k}$ is a representative of $C_{k}$. In particular,
\begin{equation}\label{2}
\widehat{K}\widehat{K^{-1}}=\sum_{k}\alpha_{k}\widehat{C}_{k} \, \, \, \, \, {\rm with} \, \, \, \, \, \alpha_{k}=\frac{|K|^2}{|G|}\sum_{\chi\in {\rm Irr}(G)}\frac{|\chi(x)|^{2}\chi(c_{k}^{-1})}{\chi(1)}.
 \end{equation} If we pour out $|\chi(x)|^{2}$ from b) we obtain
 \begin{equation}\label{3}
 |\chi(x)|^{2}=\frac{(|K|-1)\chi(1)(\chi(d)+\chi(d^{-1}))+2\chi(1)^2}{2|K|},
  \end{equation}
  and by replacing Eq.(\ref{3}) in Eq.(\ref{2}) and making easy calculations, it follows that

 \begin{equation*}
 \alpha_{k}=\frac{|K|^2}{|G|}(\sum_{\chi\in {\rm Irr}(G)} \frac{(|K|-1)(\chi(d)+\chi(d^{-1}))\chi(c_k^{-1})}{2|K|}+\sum_{\chi\in {\rm Irr}(G)} \frac{\chi(1)\chi(c_k^{-1})}{|K|}).
\end{equation*}Consequently, by using the second orthogonality relation, if $D\neq D^{-1}$, we deduce
\begin{center}
$\alpha_{k}=\left\{ \begin{array}{lc}
 			|K| &   {\rm if}  \, \, C_{k}=1 \\
 			\\
             \frac{|K|(|K|-1)}{2|D|} & \, \, \, \, \, \, \, \, \, \, \, \, \, \, \,  {\rm if}\, \, C_{k}=D\, \,  {\rm or}\, \,D^{-1} \\
             \\ 0 &   \, \, \, \, \, \, \, \, \, {\rm in} \, \, {\rm other} \, \, {\rm case}
             \end{array}
   \right.$
   \end{center} This means that $$\widehat{K}\widehat{K^{-1}}=|K|\widehat{1}+ \frac{|K|(|K|-1)}{2|D|} \widehat{D}+ \frac{|K|(|K|-1)}{2|D|} \widehat{D^{-1}}$$ and in particular, $KK^{-1}=1\cup D\cup D^{-1}$, so a) is proved. If $D=D^{-1}$, $\alpha_k$ takes the same values except  $\alpha_{k}=\frac{|K|(|K|-1)}{|D|}$ when $C_{k}=D$. In this case, the result follows by arguing similarly. $\Box$ \\

The following elementary lemma concerning commutators is basic for proving Theorem C.\\

\textbf{Lemma 2.8.} {\it Let $G$ be a finite group and let $K=x^G$, $D=d^G$ where $x$ and $d$ are elements of $G$ and $KK^{-1}=1\cup D$. Then, $\langle D\rangle=[x, G]$ and $\langle K\rangle=\langle x\rangle[x, G]$.}\\

{\it Proof.}
If $K=\lbrace x_{1}, \ldots, x_{n}\rbrace$, then $K^{-1}K=x_{1}^{-1}K\cup \cdots \cup x_{n}^{-1}K$. If $y \in x_{i}^{-1}K$, then $y=x_{i}^{-1}x_{i}^{g} \in [x_{i}, G]$ for some $g\in G$. If $i\neq j$, then $x_{j}=x_{i}^{h}$ for some $h\in G$. Thus, $[x_{j}, G]=[x_{i}^{h}, G]=[x_{i}, G]^{h}=[x_{i}, G]$. Consequently, $KK^{-1}\subseteq [x, G]$ and $\langle D\rangle\subseteq [x, G]$. On the other hand, since any element $[x,t]$ lies in $KK^{-1}$ for all $t\in G$, then $[x, G]\subseteq \langle KK^{-1}\rangle=\langle D\rangle$ and hence, $\langle D\rangle=[x, G]$. The equality $\langle K\rangle=\langle x\rangle[x, G]$ is standard, so the lemma is proved.

\section{Proof of Theorem A}

Before proving Theorem A, we analyze a particular case under the assumption $KK^{-1}=1 \cup D\cup D^{-1}$ appearing in Theorem A. If, in addition, we assume that $K=K^{-1}$, we prove in Lemma 3.1 that $D=D^{-1}$, that is, $K^2=1\cup D$, and there is no need to use the CFSG to show the non-simplicity of $G$. In fact, $\langle K\rangle$ is solvable and its structure is completely determined by the authors in Theorem A of \cite{Nuestro}.\\

\textbf{Lemma 3.1.} {\it Let $K$ and $D$ conjugacy classes of a finite group $G$ such that $KK^{-1}=1\cup D\cup D^{-1}$. If $K$ is real, then $D$ is real.}\\

{\it Proof.} Assume that $K=K^{-1}$ and let $x\in K$. If $o(x)=2$ we can assume that any two elements of $K$ do not commute, because otherwise the elements of $D$ (and of $D^{-1}$) would have order $2$ and $D$ would trivially be a real class. If $|K|=2$, say for instance $K=\{x, x^g\}$, then $KK^{-1}=1\cup \{xx^g, x^gx\}$. Notice that $\{xx^g, x^gx\}$ cannot be decomposed into two central classes. In fact, if $xx^g$ and $x^gx$ are central elements, then $x$ and $x^g$ would commute. Therefore, $|K|\geq 3$ and we write $K=\{x_1, \ldots, x_n\}$ with $n\geq3$. Let $x_{i}, x_{j}\in K$ and we have $x_{i}^{x_{j}}=x_{l}\in K$ for some positive integer $l$. Furthermore, $x_{i}\neq x_{l}$, otherwise $x_{i}$ and $x_{j}$ would commute. Thus, $x_{l}^{x_{j}}=(x_{i}^{x_{j}})^{x_{j}}=x_{i}$ and $(x_{i}x_{l})^{x_{j}}=x_{l}x_{i}=(x_{i}x_{l})^{-1}$. Since $x_{i}x_{l}\in D$ or $x_{i}x_{l}\in D^{-1}$, we conclude that $D$ is real too.\\

Suppose now that $o(x)>2$. We can clearly assume that $x^2\in D$ (analogous if $x^2\in D^{-1}$). In addition, since there exists $g\in G$ such that $x^g=x^{-1}$, we have $(x^2)^g=(x^g)^2=(x^{-1})^2=x^{-2}$. Now, if $x^2=x^{-2}$, then $o(x^{2})=2$ and $D$ is real and if $x^2\neq x^{-2}$, then $x^2$, $x^{-2}\in D$ and $D$ is real.$\Box$\\

{\it Proof of Theorem A.} Let $x, d \in G$ such that $K=x^{G}$, $D=d^G$ and $KK^{-1}=1\cup D\cup D^{-1}$. We suppose that $G$ is simple and we will look for a contradiction. We distinguish three parts appealing to the CFSG. We show that for any alternating group, simple group of Lie type or sporadic group there is no conjugacy class satisfying the hypotheses of the theorem. \\

{\it Case 1. } Suppose that $G=A_{n}$ with $n\geq 5$.\\

It is easy to check that $A_5$ does not satisfy the property of the statement for any non-trivial conjugacy class $K$. Suppose that $n>5$. Note that $x$ and $x^{-1}$ are permutations of the same type. It follows that $x^{S_{n}}(x^{-1})^{S_{n}}=1 \cup D^{S_{n}}\cup (D^{-1})^{S_{n}}=1 \cup D^{S_{n}}$ and hence $\eta(x^{S_{n}}(x^{-1})^{S_{n}})=2$. By applying Theorem 2.1, we get a contradiction because $x$ and $x^{-1}$ should be permutations of different type.\\

{\it Case 2.} Suppose that $G$ is a finite simple group of Lie type.\\

If $G$ is a finite simple group of Lie type in characteristic $p$, we can always take the Steinberg character $\psi \in $ Irr$(G)$ which satisfies $\psi(t)=\pm |{\rm \textbf{C}}_{G}(t)|_{p}$ for every $p$-regular element $t\in G$ and $\psi(t)=0$ for every $p$-singular element $t\in G$. Furthermore, $\chi(1)=|G|_p$ (see for instance Chapter 6 of \cite{Carter}). Assume that there exists a non-trivial pair of elements $x,d\in G$ such that the assertion b) of Theorem B holds and we will work to get a contradiction.\\

{\it Case 2.1.} Suppose that $x$ is $p$-regular. We know that $\psi(x)=\pm |{\rm \textbf{C}}_{G}(x)|_p\neq 0$. By the equivalence of Theorem B we have
 \begin{equation}\label{4}
|K||{\rm \textbf{C}}_{G}(x)|_{p}^{2}-|G|_{p}^{2}=\frac{|K|-1}{2}|G|_{p}(\psi(d)+\psi(d^{-1})).
\end{equation}
If $\psi(d)=\psi(d^{-1})=0$, then $|K|=|K|_{p}^{2}$ and this contradicts $p^a$-Burnside's Lemma (see for instance Theorem 15.2 of \cite{Huppert}). Thus, $\psi(d)=\psi(d^{-1})=\pm |{\rm \textbf{C}}_{G}(d)|_{p}$ and by replacing in Eq.(\ref{4}) we obtain $$(|K|_{p'}-|K|_{p})|{\rm \textbf{C}}_{G}(x)|_{p}=(|K|-1)(\pm|{\rm \textbf{C}}_{G}(d)|_{p}).$$
If $p$ divides $|K|$, it follows that $|K|_{p'}-|K|_{p}=|K|-1$, which implies $|K|=|K|_{p'}$ and $|K|_{p}=1$, a contradiciton. Consequently, $p$ does not divide $|K|$ and since $\psi(1)=|G|_{p}$ we conclude by Lemma 2.6 that either $\psi(x)=0$ or $1\neq x\in {\rm \textbf{Z}}(\psi)$. Both possibilities yield to a contradiction. \\

{\it Case 2.2.} Suppose that $x$ is $p$-singular. We know that $\psi(x)=0$ and $\psi(1)=|G|_p$. By the assertion b) of Theorem B, $$\psi(d)+\psi(d^{-1})=\frac{-2|G|_{p}}{(|K|-1)}< 0.$$

This means that $d$ is a $p$-regular element and we necessarily have
\begin{equation*}
\psi(d)=\psi(d^{-1})=-|{\rm \textbf{C}}_{G}(d)|_{p}.
\end{equation*}
As a consequence, by the two equalities above, $|K|=|D|_{p}+1$ and, thus,  $p$ does not divide $|K|$.\\

Now we prove that $x$ is a $p$-element. We consider de decomposition $x=x_{p}x_{p'}$. Notice that ${\rm \textbf{C}}_{G}(x)={\rm \textbf{C}}_{G}(x_p)\cap {\rm \textbf{C}}_{G}(x_{p'})\subseteq {\rm \textbf{C}}_{G}(x_{p'})$, which shows that $|x_{p'}^{G}|$ divides $|K|$, and then $p$ does not divide $|x_{p'}^{G}|$ either. By applying Lemma 2.6 again, we obtain that etiher $\psi(x_{p'})=0$, which leads to a contradiction because $x_{p'}$ is $p$-regular, or $x_{p'}\in {\rm \textbf{Z}}(\psi)=1$. Consequently, $x$ is a $p$-element. Since $d$ is $p$-regular, we apply Lemma 2.4 and this straightforwardly contradicts the simplicity of $G$.\\

{\it Case 3.} Suppose that $G$ is a sporadic finite simple group.\\

By using the character tables of the sporadic groups (for instance in {\sf GAP} \cite{GAP}) we can check that the equivalence of Theorem B does not hold for any of these groups and any two non-trivial conjugacy classes of it. In fact, for any sporadic simple group, the only character satisfying such assertion for fixed elements $x, d\in G$ with $x\neq1$ is the principal character.\\

The non-simplicity of $G$ when $KK^{-1}=1\cup D$ is a direct consequence of our previous arguments when $D=D^{-1}$ taking into account the corresponding case of Theorem B. $\Box$

\section{Proofs of Theorems C and D}
{\it Proof of Theorem C.} Let $K=x^G$ with $x\in G$. We write $\widehat{K}\widehat{K^{-1}}=|K|\widehat{1}+m\widehat{D}$, so $|K|^2=|K|+m|D|$ and $|D|$ divides $|K|(|K|-1)$. The fact that $\langle K\rangle/\langle D\rangle$ is cyclic follows immediately from Lemma 2.8.\\

1) Suppose that $|D|=|K|-1$. Then $|KK^{-1}|=|K|$. Note that $xK^{-1}\subseteq KK^{-1}$ and, since $|xK^{-1}|=|K|$, we obtain $xK^{-1}=KK^{-1}$. Then $K^{-1}=x^{-1}KK^{-1}$, which implies that $K^{-1}=\langle x^{-1}K\rangle K^{-1}$. This means that $K^{-1}$ is union of right classes of $\langle x^{-1}K\rangle $. Also, $\langle x^{-1}K\rangle=\langle KK^{-1}\rangle=\langle D\rangle$, so we get that $|\langle D\rangle|$ divides $|K|$. As $|K|=|KK^{-1}|\leq |\langle KK^{-1}\rangle|=|\langle D\rangle|$, then $|\langle D\rangle|=|K|$. Since $x^{-1}K\subseteq \langle KK^{-1}\rangle$ and $|x^{-1}K|=|K|=|\langle KK^{-1}\rangle|$, we obtain $\langle D\rangle=x^{-1}K$. Thus, $\langle D\rangle=xK^{-1}\subseteq 1\cup D\subseteq \langle D\rangle$, so $\langle D\rangle=1\cup D$ is $p$-elementary abelian for some prime $p$. As $\langle K\rangle/\langle D\rangle$ is cyclic, this case is finished.\\

2) Assume that $|K|=|D|$. Observe that $xK^{-1}\cup x^{-1}K\subseteq K^{-1}K$. We divide the proof of this case into two subcases: whether $xK^{-1}=x^{-1}K$ or not. Suppose first that $xK^{-1}=x^{-1}K$. We have $K=x^{2}K^{-1}$ and, analogously, $K^{-1}=(x^g)^{-2}K$ for every $g\in G$. By replacing $K^{-1}$ in the former equality, we deduce that $K=x^{2}(x^g)^{-2}K$ for every $g\in G$. We define $$N=\langle x^{2}(x^g)^{-2}\, \, | \, \, x\in K \, \, , \, \, g\in G\rangle.$$ Then $K=NK$ and, as a consequence, $|N|$ divides $|K|$. In addition, $KK^{-1}=NKK^{-1}$, so $|N|$ also divides $|KK^{-1}|=1+|D|=1+|K|$, which allows to $N=1$. As a result, $x^{2}\in {\rm \bf Z}(G)$. If $y\in K$, then $y=x^g\in K$ for some $g\in G$ and note that $y^2=(x^g)^2=(x^2)^g=x^2$. On the other hand, we can write $KK^{-1}=xK^{-1}\cup \{z\}$ for some $z\in D$. Since $xK^{-1}=x^{-1}K=(xK^{-1})^{-1}$ and $KK^{-1}$ coincides with its inverse, both facts show that $z=z^{-1}$, that is, $z$ has order 2. Now, if we take two distinct elements $y,y^g\in K$ with $g\in G$, then $y^{-1}y^g\in D$, so we write $y^g=yd$ for some $d\in D$. Then $y^2=(y^g)^2=(yd)^2=yy^d$ and consequently, $y=y^d$. This means that $[y,d]=1$ and hence $[y,y^g]=1$. Therefore, $\langle K\rangle$ is abelian, so the assertion 2) holds.\\

Assume now that $xK^{-1}\neq x^{-1}K$. We know that $xK^{-1}\cup x^{-1}K\subseteq KK^{-1}$. Since $|KK^{-1}|=|K|+1$ and $|K|=|xK^{-1}|=|x^{-1}K|$, there exists only just one element $z\in xK^{-1}\setminus x^{-1}K$. Moreover, it is easy to prove that $z^{-1}$ is the only element contained in $x^{-1}K\setminus xK^{-1}$ (notice that $z\neq z^{-1}$). Therefore, $KK^{-1}$ can be decomposed as $$KK^{-1}=xK^{-1}\cup x^{-1}K=(xK^{-1}\cap x^{-1}K)\cup \lbrace z\rbrace\cup \lbrace z^{-1}\rbrace.$$ From this equality and the fact that $(x^{-1}K)(xK^{-1})=KK^{-1}$, we deduce that
$$(KK^{-1})^2=(KK^{-1})\cup \lbrace z^{2}\rbrace\cup\lbrace z^{-2}\rbrace=1 \cup D\cup \lbrace z^{2}\rbrace\cup\lbrace z^{-2}\rbrace.$$ On the other hand, $(KK^{-1})^2=(1\cup D)(1\cup D)=1\cup D\cup D^{2}$. It follows that $D^{2}\subseteq 1 \cup D\cup \lbrace z^{2}\rbrace\cup\lbrace z^{-2}\rbrace$. We distinguish two cases. If $z^{2}\in D$, then $D^2=1\cup D$ and hence, $\langle D\rangle$ is $p$-elementary abelian for some prime $p$. We get the assertion 2) by taking into account that $\langle K\rangle/\langle D\rangle$ is cyclic. Assume now that $z^{2}\not \in D$. Then $\langle z^2\rangle\unlhd G$ because either $\{z^{2}\}$ and $\{z^{-2}\}$ are central conjugacy classes or $\{z^2, z^{-2}\}$ is a conjugacy class. We write $\overline{G}=G/\langle z^2\rangle$ and we obtain $\overline{D^{2}}\subseteq \overline{1}\cup \overline{D}$. So we have two possibilities: $\overline{D^{2}}=\overline{1}$ or $\overline{D^{2}}=\overline{1}\cup \overline{D}$. If $\overline{D^{2}}=\overline{1}$, then $\langle \overline{D}\rangle\cong \mathbb{Z}_{2}$, and as a result $\langle D\rangle$ is metacyclic. Consequently, $\langle K\rangle$ is solvable with $dl(\langle K\rangle)\leq 3$. Finally, if $\overline{D^{2}}=\overline{1}\cup \overline{D}$, it certainly follows that $\langle \overline{D}\rangle$ is elementary abelian. Then $\langle D\rangle$ is metabelian and, again $\langle K\rangle$ is solvable with $dl(\langle K\rangle)\leq 3$.\\

3) Assume that $|D|=|K|(|K|-1)$. Since $\widehat{K}\widehat{K^{-1}}=|K|\widehat{1}+m\widehat{D}$, we necessarily have $m=1$. We write $K=\{x_1,\ldots, x_{n}\}$ and $K^{-1}K=x_{1}^{-1}K\cup \ldots \cup x_{n}^{-1}K$. Notice that $1\in x_{i}^{-1}K$ for all $i=1,\cdots, n$. We rewrite the previous equality as a disjoint union
$$K^{-1}K=1\cup (x_{1}^{-1}K\setminus 1)\cup \ldots \cup (x_{n}^{-1}K\setminus 1).$$ By counting elements we conclude that $x_{i}^{-1}K\cap x_{j}^{-1}K=1$ for all $i=1,\ldots, n$ with $i\neq j$. Let $g\in {\rm \bf C}_{G}(x_{i}x_{j}^{-1})$ with $i\neq j$. Thus, $(x_{i}x_{j}^{-1})^g=x_{i}^{g}(x_{j}^{-1})^{g}=x_{i}x_{j}^{-1}$. From the last equality we have $x_{i}^{-1}x_{i}^{g}=x_{j}^{-1}x_{j}^{g}=1$, so $g\in {\rm \bf C}_{G}(x_{i})\cap {\rm \bf C}_{G}(x_{j}^{-1})$. Hence, ${\rm \bf C}_{G}(x_{i}x_{j}^{-1})={\rm \bf C}_{G}(x_{i})\cap {\rm \bf C}_{G}(x_{j}^{-1})$. As $x_{i}x_{j}^{-1}\in {\rm \bf C}_{G}(x_{i}x_{j}^{-1})$, then $x_{i}x_{j}^{-1}\in {\rm \bf C}_{G}(x_{i})$, so $[x_{i}, x_{j}]=1$. Therefore, $\langle K\rangle$ is generated by pairwise commuting elements, which means that $\langle K\rangle$ is abelian.$\Box$\\

{\it Proof of Theorem D.} We write
\begin{equation}\label{A}
\widehat{K}\widehat{K^{-1}}=|K|\widehat{1}+m\widehat{K}+m\widehat{K^{-1}}
\end{equation}
where $m$ is a positive integer. Note that in the product $\widehat{K}\widehat{K^{-1}}$, the sums $\widehat{K}$ and $\widehat{K^{-1}}$ necessarily have the same multiplicity. Let us calculate $\widehat{K}^2$ and $(\widehat{K^{-1}})^2$. By Lemma 2.7(3), we know that $(\widehat{K}^2, \widehat{K})=(\widehat{K}\widehat{K^{-1}}, \widehat{K})=m=(\widehat{K}\widehat{K^{-1}}, \widehat{K^{-1}})=(\widehat{K}^2, \widehat{K^{-1}})$, so
\begin{equation}\label{A1}
\widehat{K}^2=m\widehat{K}+m\widehat{K^{-1}}+\widehat{T_1}
\end{equation}
where $T_{1}$ is union of classes, repeated or not, different from $K$ and $K^{-1}$, that is, $(\widehat{T_{1}}, \widehat{L})=0$ if $L\in \lbrace K, K^{-1}\rbrace$. Analogously, $$(\widehat{K^{-1}})^2=m\widehat{K}+m\widehat{K^{-1}}+\widehat{T_2}$$ with  $(\widehat{T_{2}}, \widehat{L})=0$ if $L\in \lbrace K, K^{-1}\rbrace$.\\

Now, let us see that if $\widehat{T}_{1}=0$, then the theorem is proved. If $\widehat{T}_{1}=0$, then $K^{3}=KK^2=K(K\cup K^{-1})=1\cup K\cup K^{-1}=K^{n}$ for every $n\geq 2$, so $\langle K \rangle=1\cup K \cup K^{-1}$. In particular, $\langle K\rangle$ is $p$-elementary abelian for some prime $p$ and we write $|\langle K\rangle|=p^n$. Notice that $p\neq 2$ because otherwise $K$ would be real against the hypotheses. Furthermore, by Eq.(\ref{A}), $|K|^{2}=|K|+2m|K|$ and this forces $|K|$ to be odd. Also, $p^{n}=|\langle K\rangle|=1+2|K|$ and hence  $$|K|=\frac{p-1}{2}(p^{n-1}+p^{n-2}+\cdots + 1),$$ so both factors of this product should be odd. In particular, $\frac{p-1}{2}=2r+1$ with $r$ a positive integer, that is, $p\equiv3$ mod $4$. In addition, $n$ is odd as the second factor of the product consists of $n$ odd addends. Analogously, the theorem can be proved when $\widehat{T}_{2}=0$.\\

Therefore, from now on, we may assume that $\widehat{T}_{1}\neq 0\neq \widehat{T}_{2}$ and we see that this leads to a contradiction. Let us show first that $\widehat{T_{1}}=\widehat{T_{2}}$. We have $$\widehat{K}(\widehat{K}\widehat{K^{-1}})=\widehat{K}(|K|\widehat{1}+m\widehat{K}+m\widehat{K}^{-1})=|K|\widehat{K}+2m^2\widehat{K}+2m^{2}\widehat{K}^{-1}+m|K|\widehat{1}+m\widehat{T_{1}},$$
and on the other hand
$$\widehat{K}^2\widehat{K^{-1}}=(m\widehat{K}+m\widehat{K^{-1}}+\widehat{T_1})\widehat{K}^{-1}=m|K|\widehat{1}+2m^2\widehat{K}+2m^2\widehat{K}^{-1}+m\widehat{T_{2}}+\widehat{T_{1}}\widehat{K}^{-1}.$$
Comparing both equalities we deduce that
\begin{equation} \label{B}
\widehat{T_{1}}\widehat{K^{-1}}=|K|\widehat{K}+m\widehat{T_{1}}-m\widehat{T_{2}}.
\end{equation} By reasoning similarly with $\widehat{K^{-1}}(\widehat{K}\widehat{K^{-1}})$ and $(\widehat{K^{-1}})^2\widehat{K}$ we get
\begin{equation}\label{C}
\widehat{T_{2}}\widehat{K}=|K|\widehat{K^{-1}}+m\widehat{T_{2}}-m\widehat{T_{1}}.
\end{equation}
 By adding both equalities,
\begin{equation*}
\widehat{T_{1}}\widehat{K^{-1}}+\widehat{T_{2}}\widehat{K}=|K|\widehat{K}+|K|\widehat{K^{-1}},
\end{equation*}
which implies that $\widehat{T_{1}}\widehat{K^{-1}}=m_{1}\widehat{K}+m_{2}\widehat{K^{-1}}$, where $m_{1}$ and $m_{2}$ are positive integers. By comparing this equation with Eq.(\ref{B}), it necessarily follows that $m_{1}=|K|$, $m_{2}=0$ and $\widehat{T}_{1}=\widehat{T}_{2}$ as we wanted to prove. We conclude that $\widehat{T_{1}}\widehat{K^{-1}}=|K|\widehat{K}$. Similarly we can deduce that $\widehat{T_{2}}\widehat{K}=|K|\widehat{K^{-1}}=\widehat{T_{1}}\widehat{K}$. As a result, if we denote $T=T_1=T_2$, we have
$$TK^{-1}=K$$
$$TK=K^{-1}.$$
On the other hand, by Eq.(\ref{A1}),
\begin{equation}\label{D}
K^{2}=K\cup K^{-1}\cup T,
\end{equation}
so $|K^2|=2|K|+|T|$. In addition, $TK^{2}=(TK)K=K^{-1}K=1\cup K\cup K^{-1}$, which implies that $1+2|K|=|TK^{2}|\geq |K^2|=2|K|+|T|$. As a consequence, $|T|=1$ and we write $T=\lbrace t\rbrace$. Observe that $t\neq 1$ because otherwise, by Eq.(\ref{D}), $K$ would be real. Then $\langle K\rangle=1\cup K \cup K^{-1}\cup \lbrace t\rbrace$, so $|\langle K\rangle\cap {\rm \textbf{Z}}(G)|=2$ and $o(t)=2$. Therefore, the non-trivial elements of $\langle K\rangle$ have order 2 or $s$ and this forces that $s=4$, otherwise there would appear more possible orders for the elements of $\langle K\rangle$. Since $\langle K\rangle$ is a $2$-group with only one involution, by applying Lemma 2.5, either $\langle K\rangle\cong \mathbb{Z}_{4}$ or $\langle K\rangle\cong Q_{8}$, the quaternion group of order $8$. In the first case, $|K|=1$ and this is a contradiction. If $\langle K\rangle\cong Q_{8}$, then $8=2+2|K|$, so $|K|=3$. However, we know that $|x^{\langle K\rangle}|=2$ for every $x\in K\cup K^{-1}$. Since $|x^{\langle K\rangle}|$ divides $|x^G|=|K|$, we get a contradiction, as wanted. $\Box$

\section{Examples}

We emphasize the feasibility of our results by giving examples. Some of them have been found by using the {\sc SmallGroups}  library of {\sf GAP} \cite{GAP}. The $m$-th group of order $n$ in this library is identified by $n\#m$. First, we provide an example illustrating the non-simplicity of a group satisfying the hypothesis of Theorem A.\\

{\bf Example 5.1.} Let $G=\langle a\rangle\rtimes \langle b\rangle$ where $\langle a\rangle\cong \mathbb{Z}_{7}$, $\langle b\rangle\cong \mathbb{Z}_{3}$ and $a^b=a^2$. Let $K$ be one of the two classes of elements of order 3, which satisfies $|K|=7$. It holds $KK^{-1}=1\cup D\cup D^{-1}$ where $D$ is a conjugacy class of elements of order 7 and size 3. We have $\langle K\rangle=G$ and $\langle D\rangle\cong  \mathbb{Z}_{7}$. This is the example of a group with the smallest order satisfying the property asserted in Theorem A with $D\neq D^{-1}$.\\

The three following examples correspond to Theorem C.\\

{\bf Example 5.2.} Let $G=(\langle c\rangle \rtimes \langle b\rangle)\rtimes \langle a\rangle$ such that $\langle a\rangle\cong \mathbb{Z}_{2}$, $\langle b\rangle\cong \mathbb{Z}_{5}$ and $\langle c\rangle\cong \mathbb{Z}_{11}$ with the following actions: $c^a=c^{10}$, $b^a=b$ and $c^b=c^4$. This is the group $110\#1$ of {\sf GAP}. Take $K=b^G$ which satisfies $KK^{-1}=1\cup D$, $|K|=11$ and $|D|=10$ where $D=c^G$. We have $\langle K\rangle=\langle c\rangle \rtimes \langle b\rangle\cong {\mathbb Z}_{11}\rtimes {\mathbb Z}_{5}$ and $\langle D\rangle=\langle c\rangle\cong {\mathbb Z}_{11}$. This is an example of Case 1 of Theorem C.\\

{\bf Example 5.3.} Let $G=\langle a\rangle \times A_{4}$ where $\langle a\rangle \cong {\mathbb Z}_{10}$ and we take $K=x^{G}$ where $x=at$ and $t$ is an involution of $A_4$. Thus, $o(x)=10$, $|K|=3$ and $KK^{-1}=1\cup D$ where $D=t^{G}$ and $|D|=3$. In this example, $\langle K\rangle\cong \mathbb{Z}_{10}\times \mathbb{Z}_{2}\times \mathbb{Z}_{2}$ and $\langle D\rangle\cong \mathbb{Z}_{2}\times \mathbb{Z}_{2}$. This is an example of Case 2 of Theorem C.\\

{\bf Example 5.4.} Let $G=((\langle c\rangle \times \langle d\rangle)\rtimes \langle b\rangle)\rtimes \langle a\rangle$ such that $\langle a\rangle\cong \mathbb{Z}_{2}$, $\langle b\rangle\cong \mathbb{Z}_{3}$ and $\langle c\rangle\cong \mathbb{Z}_{5}\cong \langle d \rangle$ with the following actions: $c^b=cd^3$, $d^b=c^4d^3$, $c^a=c$, $d^a=c^4d^4$ and $b^a=b^2$. This is the group $150\#5$ of {\sf GAP}. Take $K=c^G=\{c, cd^3, c^3d^2\}$ which satisfies $KK^{-1}=1\cup D$, $|K|=3$ and $|D|=6$ where $D=d^G=\{d,d^4, c^4d^4,cd,c^4d^3,cd^2\}$. We have $\langle K\rangle=\langle c\rangle\times \langle d\rangle\cong {\mathbb Z}_{5}\times {\mathbb Z}_{5}$. This is an example of Case 3 of Theorem C.\\

In the next example, $|D|$ is a divisor of $|K|(|K|-1)$ different from those appearing in Theorem C.\\

{\bf Example 5.5.} Let $G=SL(2,3)$. Take $K=x^G$ where $o(x)=3$. It is satisfied $|K|=4$ and $KK^{-1}=1\cup D$, where $D$ is a class of elements of order 4 and $|D|=6$. We have $\langle K\rangle \cong SL(2,3)$ and $\langle D\rangle\cong Q_{8}$. Note that $dl(\langle K\rangle)=3$.\\

We give two examples illustrating Theorem D. In the first one, $\langle K\rangle$ is cyclic whereas in the second it is not.\\

{\bf Example 5.6.} Let $G=\mathbb{Z}_{11}\rtimes\mathbb{Z}_{5}$ the semidirect product of $\langle a\rangle\cong \mathbb{Z}_{5}$ acting non-trivially on $\langle x\rangle\cong \mathbb{Z}_{11}$ by $x^a=x^4$. Take the conjugacy class $K=x^G=\lbrace x, x^3, x^4, x^5, x^9\rbrace$. It follows that $KK^{-1}=1 \cup K\cup K^{-1}=\langle x\rangle$ with $|\langle K\rangle|=11$ and $|K|=5$.\\

{\bf Example 5.7.} Let $N=\langle x\rangle \times \langle y\rangle \times \langle z\rangle\cong \mathbb{Z}_{3}\times  \mathbb{Z}_{3}\times  \mathbb{Z}_{3}$ and $H=\langle a\rangle\cong  \mathbb{Z}_{13}$. Let us consider $a\in {\rm Aut}(N)$ defined as follows: $x^a=y$, $y^a=z$ and $z^a=xy$. Write $G=N\rtimes H$ the semidirect product of $N$ and $H$. This is the group $351\#12$ of {\sf GAP}. Take $$K=x^G=\lbrace x, y, z, xy, yz, xyz, xy^2z, xy^2z^2, x^2z^2, x^2y, y^2z, xyz^2, x^2z\rbrace.$$ We have $|\langle K\rangle|=27$ and $|K|=13$.

\section{Analogous problems for irreducible characters}
Following the parallelism between conjugacy classes and irreducible characters, we reflect on the problem of translating our results into Character Theory. As we have asserted in the Introduction by means of an example, the fact that $\eta(KK^{-1})=3$ does not imply the non-simplicity of the group. Something similar occurs when working with irreducible characters. For example, if we consider the simple group PSL$(2,11)$, there exist three irreducible characters $\chi$, $\psi$ and $\varphi$ such that $\chi\overline{\chi}=1+\psi+\varphi$ (see for instance page 290 of \cite{Isaacs}).\\

Trying to transfer Theorem A into the framework of irreducible characters, we find that in \cite{AdanBante2} the author gives the structure of a finite solvable group $G$ with $\chi \in$ Irr$(G)$ such that $\chi\overline{\chi}=1_G+m_1\alpha_1+m_2\alpha_2$ where $\alpha_1$, $\alpha_2\in$ Irr$(G)$ are non-principal characters and $m_1$ and $m_2$ are strictly positive integers. We are not aware, however, whether the above equality may hold in a simple group when $\alpha_{2}=\overline{\alpha_{1}}$. On the other hand, regarding the particular case of Theorem A, if we take $G={\rm PSp}_{2n}(3)$, $n\geq 2$, or $G={\rm PSU}_{\rm n}(2)$, $(n,3)=1$, $n\geq 4$, it is known that there exists a non-trivial character $\psi \in {\rm Irr}(G)$ such that $\psi\overline{\psi}=\chi+1$, with $\chi \in {\rm Irr}(G)$ (see \cite{Malle} for instance). So we conclude that simplicity may occur when the particular case of Theorem A is translated into irreducible characters. \\

The statement of Theorem D has no sense when dealing with irreducible characters. In fact, let $\chi \in {\rm Irr}(G)$ such that $\chi \overline{\chi}=a1_{G}+b\chi+c\overline{\chi}$ with $a, b, c \in \mathbb{Z}^{+}$. Then $a=[\chi\overline{\chi}, 1_{G}]=[\chi, \chi]=1$ and $b=[\chi\overline{\chi}, \chi]=[\chi\overline{\chi}, \overline{\chi}]=c$, so $\chi\overline{\chi}=1_{G}+b\chi+b\overline{\chi}$. By taking degrees, $\chi(1)^{2}=1+2b\chi(1)$, so $\chi(1)(\chi(1)-2b)=1$. This forces that $\chi(1)=1$ and $\chi(1)=2b+1$, a contradiction.\\

\noindent {\bf Acknowledgements}

\bigskip
The results in this paper are part of the third author's Ph.D. thesis, and she acknowledges the predoctoral grant PREDOC/2015/46, Universitat Jaume I. The first and second authors are supported by the Valencian Government, Proyecto PROMETEOII/2015/011. The first and the third authors are also partially supported by Universitat Jaume I, grant P11B\-2015-77.

\end{document}